\documentclass[12pt]{article}

\usepackage{amsmath}

\title{On the derivatives of the Heun functions}
\author{Galina Filipuk\footnote{Faculty of Mathematics, Informatics and Mechanics,
University of Warsaw, Banacha 2, Warsaw, 02-097, Poland.
Email: G.Filipuk@mimuw.edu.pl}, Artur Ishkhanyan\footnote{Russian-Armenian University, Yerevan, 0051 Armenia}
\footnote{Institute for Physical Research of NAS of Armenia, 0203 Ashtarak, Armenia. Email: aishkhanyan@gmail.com}, Jan Derezi\'nski\footnote{Department of Mathematical Methods in Physics, Faculty of Physics, University of Warsaw, Pasteura 5, 02-093, Warszawa, Poland. Email: jan.derezinski@fuw.edu.pl}}
\date{\today}

\begin{document}
\maketitle

\begin{abstract}
The Heun functions satisfy linear ordinary differential equations of second order with certain singularities in the complex plane.
The first order derivatives of the Heun functions satisfy linear second order differential equations with one more singularity. In this paper we  compare  these equations with  linear differential equations  isomonodromy deformations of which are described  by
 the Painlev\'e equations $P_{II}-P_{VI}$.
\end{abstract}



\section{Introduction}

The general Heun equation is the most general second-order linear  Fuchsian  ordinary differential equation with four regular singular points in the complex plane \cite{Heun, Ronveaux, Sl_Lay, dlmf}.  Although it  is a genaralization of the well-studied Gauss hypergeometric equation with three regular singularities, it is much more difficult to investigate   properties of the Heun functions. The additional singularity  causes many complications in comparison with the hypergeometric case (for instance, the solutions in general have no integral representations involving simpler mathematical functions). There also exist confluent Heun equations (see \cite{Ronveaux, Sl_Lay}) which have irregular singularities.
There are many studies on  the properties of solutions of the Heun equations from different perspectives (see, for instance, \cite{GF, Ishkhanyan, Ishkhanyan2, Ishkhanyan3, Ku, LeroyIshkhanyan, Maier, Ronveaux2, conn, RV_GF, RV_GF2} and the references therein).  The Heun functions (and their confluent cases) appear extensively in many problems of mathematics, mathematical physics, physics and engineering (e.g., \cite{App1, App2, App3}).

The general Heun equation is given by the following equation:
\begin{equation}\label{GHE}
\frac{d^2 u}{dz^2}+\left(\frac{\gamma}{z}+\frac{\delta}{z-1}+\frac{\varepsilon}{z-t}\right)\frac{du}{dz}+\frac{\alpha\beta z-q}{z(z-1)(z-t)}u=0,
\end{equation}
where the parameters satisfy the Fuchsian relation
\begin{equation}\label{GHE param}
1+\alpha+\beta=\gamma+\delta+\varepsilon.
\end{equation}
This equation has four regular singular points at $z=0,\,1,\,t$ and $\infty$. Its solutions, the Heun functions, are usually denoted by $u=H(t,\,q;\,\alpha,\,\beta,\,\gamma,\,\delta;\,z)$ assuming that $\varepsilon$ is obtained from (\ref{GHE param}). The parameter $q$ is referred to as the accessory parameter.

It is well-known that the derivative of the hypergeometric function ${}_2F_1$ is again a hypergeometric function with different values of the parameters. However, for the Heun function it is generally not the case. The first order derivative of the general Heun function satisfies a second order Fuchsian differential equation with five regular singular points \cite{Ishkhanyan, Ishkhanyan2, LeroyIshkhanyan}. It can be verified by direct computations that the function $v(z)=du/dz$, where $u=u(z)$ is a solution of (\ref{GHE}), satisfies the following equation:
\begin{equation}\label{dGHE}
\frac{d^2v}{dz^2}+\left(\frac{\gamma+1}{z}+\frac{\delta+1}{z-1}+\frac{\varepsilon+1}{z-t}-\frac{\alpha\beta}{\alpha\beta z-q}\right)\frac{dv}{dz}+
\frac{f(z)}{z(z-1)(z-t)(\alpha\beta z-q)}v=0,
\end{equation}
where $f(z)=z(\alpha\beta z-2q)(\alpha\beta+\gamma+\delta+\varepsilon)+(q^2+q(\gamma+t (\gamma+\delta)+\varepsilon)-\alpha\beta \gamma t)$. We see from that an additional singularity at $z=q/(\alpha\beta)$ involving the accessory parameter is added.

It is known that in some cases  equation (\ref{dGHE}) reduces to a Heun equation (\ref{GHE}) with altered parameters   \cite{Ishkhanyan2}. Indeed, we can  observe that in four cases when $q=0,\,\alpha\beta,\,\alpha\beta t$ and $\alpha\beta=0$ the additional singularity in (\ref{dGHE}) disappears and we obtain the Heun equation (\ref{GHE})  with different parameters \cite{Ishkhanyan2}. The equation for the  derivative of the Heun functions allows one to construct several new expansions  of solutions of the Heun equations in terms of various special functions (e.g., hypergeometric functions)    \cite{Ishkhanyan}. Similar results hold for confluent cases \cite{LeroyIshkhanyan}.

This paper is organized as follows. In Section~\ref{section: list} we give a list of all confluent Heun equations together with linear second order equations for the derivatives of the Heun functions. In Section~\ref{section: Painleve}  we briefly describe the theory of isomonodromy deformations of linear equations and show how the famous Painlev\'e equations appear in this context. Next, in Section~\ref{section: main results} we present our main results. In particular, we will compare  linear equations for the Heun  derivatives with linear  differential equations isomonodromy deformations of which are described by the Painlev\'e equations.

\section{Confluent Heun equations and equations for derivatives of confluent Heun functions}\label{section: list}

The general Heun equation is given by (\ref{GHE}) together with (\ref{GHE param}) and the linear equation for the derivative of the Heun functions is (\ref{dGHE}).

The confluent Heun equation is written as
\begin{equation}\label{CHE}
\frac{d^2u}{dz^2}+\left(\frac{\gamma}{z}+\frac{\delta}{z-1}+\varepsilon\right)\frac{du}{dz}+\frac{\alpha z-q}{z(z-1)}u=0
\end{equation}
and the linear  equation for the function $v=du/dz$ is given by
\begin{equation}\label{dCHE}
\frac{d^2v}{dz^2}+\left(\frac{\gamma+1}{z}+\frac{\delta+1}{z-1}+\varepsilon-\frac{\alpha}{\alpha z-q} \right)\frac{dv}{dz}+\frac{g(z)}{z(z-1)(\alpha z-q)}v=0,
\end{equation}
where $g(z)=(\alpha+\varepsilon)(\alpha z^2-2q z)+(q^2-(\gamma+\delta-\varepsilon)q+\alpha\gamma)$.

The double-confluent Heun equation is
\begin{equation}\label{DHE0}
\frac{d^2u}{dz^2}+\left(\frac{\gamma}{z^2}+\frac{\delta}{z}+\varepsilon\right)\frac{du}{dz}+\frac{\alpha z-q}{z^2}u=0
\end{equation}
and the linear  equation for the function $v=du/dz$ is given by
\begin{equation}\label{dDHE0}
\frac{d^2v}{dz^2}+\left(\frac{\gamma}{z^2}+\frac{\delta+2}{z}+\varepsilon-\frac{\alpha}{\alpha z-q} \right)\frac{dv}{dz}+\frac{h(z)}{z^2  (\alpha z-q)}v=0,
\end{equation}
where $h(z)=(\alpha+\varepsilon)(\alpha z^2-2q z)+(q^2-\delta q-\alpha\gamma)$.

The bi-confluent Heun equation is
\begin{equation}\label{BHE1}
\frac{d^2u}{dz^2}+\left(\frac{\gamma}{z}+ \delta +\varepsilon z\right)\frac{du}{dz}+\frac{\alpha z-q}{z}u=0
\end{equation}
and the linear  equation for the function $v=du/dz$ is given by
\begin{equation}\label{dBHE1}
\frac{d^2v}{dz^2}+\left(\frac{\gamma+1}{z}+\delta+\varepsilon z-\frac{\alpha}{\alpha z-q}\right)\frac{dv}{dz}+\frac{k(z)}{z(\alpha z-q)}v=0,
\end{equation}
where $k(z)=(\alpha+\varepsilon)z(\alpha z-2q)+(q^2-\delta q-\alpha \gamma)$.

The tri-confluent Heun equation is
\begin{equation}\label{dTHE}
\frac{d^2u}{dz^2}+\left( \gamma+\delta z+\varepsilon z^2\right)\frac{du}{dz}+(\alpha z-q)u=0
\end{equation}
and the linear  equation for the function $v=du/dz$ is given by
\begin{equation}\label{dTHE}
\frac{d^2v}{dz^2}+\left(\gamma+\delta z+\varepsilon z^2 -\frac{\alpha}{\alpha z-q}\right)\frac{dv}{dz}+\frac{p(z)}{(\alpha z-q)}v=0,
\end{equation}
where $p(z)=(\alpha+\varepsilon)(\alpha z^2-2q z)+(q^2-\delta q-\alpha \gamma)$.

\section{Isomonodromic deformations of linear equations and the Painlev\'e equations}\label{section: Painleve}

In this section we briefly review the theory of isomonodromic deformations of linear second order  differential equations following \cite{Iwasaki, Ohyama, Okamoto}. We shall use the notation similar to \cite{Ohyama}.

The isomonodromic deformations  of linear second order differential equations of the form
\begin{equation}\label{eq}
\frac{d^2v}{dz^2}+p_1(z)\frac{dv}{dz}+p_2(z)v=0
\end{equation}
with $p_1,\,p_2$ being rational functions of $z$ and parameters of deformation $t_1,\ldots,t_n$, are governed by a  completely integrable Hamiltonian system of partial differential equations with respect to the parameters. When there is one parameter of deformation, $t$,  the Painlev\'e equations $P_I-P_{VI}$ appear as the compatibility condition of the extended linear system consisting of equation  (\ref{eq}) and equaton
\begin{equation}\label{eq1}
\frac{\partial v}{\partial t}=a(z,t)\frac{\partial v}{\partial z}+b(z,t)v.
\end{equation}

The Painlev\'e equations $P_I-P_{VI}$ are nonlinear second order differential equations with the so-called Painlev\'e property. They have many interesting properties and appear in many areas of mathematics. See, for instance,  \cite{Clarkson, Gromak, Iwasaki} and numerous references therein.  The completely integrable Hamiltonian system is then equivalent to the Painlev\'e equations for one of the variables. Below we shall present necessary formulas for equations $P_{II}-P_{VI}$.

To get the sixth Painlev\'e equation one chooses
\begin{eqnarray}\label{P6 coeff p1}
p_1(z,t)&=&\frac{1-\kappa_0}{z}+\frac{1-\kappa_1}{z-1}+\frac{1-\theta}{z-t}-\frac{1}{z-\lambda},\\\label{P6 coeff p2}
p_2(z,t)&=&\frac{\kappa}{z(z-1)}-\frac{t(t-1)H_{VI}}{z(z-1)(z-t)}+\frac{\lambda(\lambda-1)\mu}{z(z-1)(z-\lambda)},
\end{eqnarray}
where
\begin{eqnarray*}
t(t-1)H_{VI}&=&\lambda(\lambda-1)(\lambda-t)\mu^2\\ &&
-\{\kappa_0(\lambda-1)(\lambda-t)+\kappa_1\lambda(\lambda-t)+(\theta-1)\lambda(\lambda-1)\}\mu+\kappa(\lambda-t).
\end{eqnarray*}
Then the compatibility between (\ref{eq})  and  (\ref{eq1}) with certain $a(z,t)$ and $b(z,t)$ (see \cite{Iwasaki, Ohyama, Okamoto} for details) leads to the
Hamiltonian system
\begin{equation*}
\frac{d\lambda}{dt}=\frac{\partial H_{VI}}{\partial \mu},\;\;\;\frac{d\mu}{dt}=-\frac{\partial H_{VI}}{\partial \lambda}
\end{equation*}
and by eliminating the function $\mu$ one can get the sixth  Painlev\'e equation
\begin{eqnarray}\label{P6}\nonumber
\frac{d^2\lambda}{dt^2}&=&\frac{1}{2}\left(\frac{1}{\lambda}+\frac{1}{\lambda-1}+\frac{1}{\lambda-t}\right)\left(\frac{d\lambda}{dt}\right)^2-\left(\frac{1}{t}+\frac{1}{t-1}+\frac{1}{\lambda-t}\right)\frac{d \lambda}{dt}
\\
&& +\frac{\lambda(\lambda-1)(\lambda-t)}{t^2(t-1)^2}\left(\alpha_6+\beta_6\frac{t}{\lambda^2}+\gamma_6\frac{t-1}{(\lambda-1)^2}+\delta_6\frac{t(t-1)}{(\lambda-t)^2}\right),
\end{eqnarray}
where
$$\alpha_6=\frac{1}{2}\kappa_{\infty}^2,\;\;\beta_6=-\frac{1}{2}\kappa_0^2,\;\;\gamma_6=\frac{1}{2}\kappa_1^2,\;\;\delta_6=\frac{1}{2}(1-\theta^2)$$
and $$\kappa=\frac{1}{4}(\kappa_0+\kappa_1+\theta-1)^2-\frac{1}{4}\kappa_{\infty}^2.$$

To get the fifth Painlev\'e equation one chooses
\begin{eqnarray}\label{P5 coeff p1}
p_1(z,t)&=&\frac{1-\kappa_0}{z}+\frac{\eta t}{(z-1)^2}+\frac{1-\theta}{z-1}-\frac{1}{z-\lambda},
\\\label{P5 coeff p2}
p_2(z,t)&=&\frac{\kappa}{z(z-1)}-\frac{t H_V}{z(z-1)^2}+\frac{\lambda(\lambda-1)\mu}{z(z-1)(z-\lambda)},
\end{eqnarray}
where
\begin{equation*}
tH_{V}=  \lambda(\lambda-1)^2\mu^2-\{\kappa_0(\lambda-1)^2+\theta \lambda(\lambda-1)-\eta t \lambda\}\mu+\kappa(\lambda-1).
\end{equation*}
Then similarly to the previous case the  corresponding Hamiltonian system with the Hamiltonian $H_V$
leads to the fifth  Painlev\'e equation
\begin{eqnarray}\label{P5}\nonumber
\frac{d^2\lambda}{dt^2}&=& \left(\frac{1}{2\lambda}+
\frac{1}{\lambda-1}\right)\left(\frac{d\lambda}{dt}\right)^2-\frac{1}{t} \frac{d \lambda}{dt}+\frac{(\lambda-1)^2}{t^2}\left(\alpha_5 \lambda+\frac{\beta_5}{\lambda}\right)\\&&
+\gamma_5\frac{\lambda}{t}+\delta_5\frac{\lambda(\lambda+1)}{\lambda-1},
\end{eqnarray}
where
$$\alpha_5=\frac{1}{2}\kappa_{\infty}^2,\;\;\beta_5=-\frac{1}{2}\kappa_0^2,\;\;\gamma_5=(1+\theta)\eta,\;\;\delta_5=\frac{1}{2}\eta^2$$
and $$\kappa=\frac{1}{4}(\kappa_0+\theta)^2-\frac{1}{4}\kappa_{\infty}^2.$$

To get the fourth Painlev\'e equation one chooses
\begin{eqnarray}\label{P4 coeff p1}
p_1(z,t)&=& \frac{1-\kappa_0}{z}-\frac{z+2t}{2}-\frac{1}{z-\lambda},
\\\label{P4 coeff p2}
p_2(z,t)&=& \frac{1}{2}\theta_{\infty}-\frac{H_{IV}}{2z}+\frac{\lambda \mu}{z(z-\lambda)},
\end{eqnarray}
where
\begin{equation*}
H_{IV}=  2\lambda\mu^2-(\lambda^2+2t \lambda+2\kappa_0)\mu+\theta_{\infty}\lambda.
\end{equation*}
Then the corresponding Hamiltonian system with the Hamiltonian $H_{IV}$
leads to the fourth  Painlev\'e equation
\begin{eqnarray}\label{P4}
\frac{d^2\lambda}{dt^2}=  \frac{1}{2\lambda} \left(\frac{d\lambda}{dt}\right)^2+\frac{3}{2}\lambda^3+4t \lambda^2+2(t^2-\alpha_4)\lambda+\frac{\beta_4}{\lambda},
\end{eqnarray}
where
$$\alpha_4=-\kappa_0+2\theta_{\infty}+1 ,\;\;\beta_4=-2\kappa_0^2.$$

The standard third Painlev\'e equation is given by
\begin{equation}\label{P3}
\frac{d^2\lambda}{dt^2}=\frac{1}{\lambda}\left(\frac{d\lambda}{dt}\right)^2-\frac{1}{t}\frac{d\lambda}{dt}+
\frac{\alpha_3\lambda^2+\beta_3}{t}+\gamma_3 \lambda^3+\frac{\delta_3}{\lambda}.
\end{equation}
However, for our purpose it is more convenient to consider equation which can be obtained from (\ref{P3}) by changing $\lambda(t)\to\lambda(t^2)/t$ and by renaming the new variable $\tau=t^2$ as $t$ again. This equation is given by
\begin{equation}\label{P3p}
\frac{d^2\lambda}{dt^2}=\frac{1}{\lambda}\left(\frac{d\lambda}{dt}\right)^2-\frac{1}{t}\frac{d\lambda}{dt}+
\frac{\alpha_3\lambda^2+\gamma_3\lambda^3}{4t^2}+\frac{\beta_3}{4t}+\frac{\delta_3}{4\lambda}.
\end{equation}
Equation (\ref{P3p}), which will be denoted by $P^{\prime}_{III},$ appears in the result of isomonodromic deformations of linear equation (\ref{eq}) with
\begin{eqnarray}\label{P3 coeff p1}
p_1(z,t)&=& \frac{\eta_0 t}{z^2}+\frac{1-\theta_0}{z}-\eta_{\infty}-\frac{1}{z-\lambda} ,
\\\label{P3 coeff p2}
p_2(z,t)&=& \frac{\eta_{\infty}(\theta_0+\theta_{\infty})}{2z}-\frac{t H^{\prime}_{III}}{z^2}+\frac{\lambda\mu}{z(z-\lambda)},
\end{eqnarray}
where
\begin{equation*}
t H^{\prime}_{III}=   \lambda^2\mu^2-\{\eta_{\infty}\lambda^2+\theta_0 \lambda-\eta_0 t\}\mu+\frac{1}{2}\eta_{\infty}(\theta_0+\theta_{\infty})\lambda
\end{equation*}
and the parameters are related by
$$\alpha_3=-4\eta_{\infty}\theta_{\infty},\;\;\beta_3=4\eta_0(1+\theta_0),\;\; \gamma_3=4\eta_{\infty}^2,\,\, \delta_3=-4\eta_0^2.$$

Finally, the second Painlev\'e equation
\begin{equation}\label{P2}
\frac{d^2\lambda}{dt^2}= 2\lambda^3+t \lambda+\alpha_2
\end{equation}
 appears in the result of isomonodromic deformations of linear equation (\ref{eq}) with
\begin{eqnarray}\label{P2 coeff p1}
p_1(z,t)&=& -2z^2-t-\frac{1}{z-\lambda},
\\\label{P2 coeff p2}
p_2(z,t)&=&  -(2\alpha_2+1)z-2H_{II}+\frac{\mu}{z-\lambda},
\end{eqnarray}
where
\begin{equation}\label{H2}
H_{II}=\frac{1}{2}\mu^2-\left(\lambda^2+\frac{1}{t}\right)\mu-\left(\alpha_2+\frac{1}{2}\right)\lambda.
\end{equation}

\section{Main results}\label{section: main results}

In this section we compare equations for the derivatives of the Heun functions with linear differential equations isomonodromy deformations of which are governed by   the Painlev\'e equations $P_{II}-P_{VI}$.

Let us consider equation for the derivative of the general Heun function  (\ref{dGHE}). By choosing parameters
\begin{gather*}
\alpha\beta=\kappa_0+\kappa_1+\theta+\kappa,\;\;\beta=\frac{1}{2}(\pm \kappa_{\infty}-1-\kappa_0-\kappa_1-\theta),\\
\gamma=-\kappa_0,\;\;\delta=-\kappa_1,\;\;\varepsilon=-\theta,\;\;q=\alpha\beta \lambda,
\end{gather*}
we can calculate that the resulting equation is the same as equation (\ref{eq}) with (\ref{P6 coeff p1}), (\ref{P6 coeff p2}) and the expression for $H_{VI}$  provided that
$$\mu=\frac{\kappa_0}{\lambda}+\frac{\kappa_1}{\lambda-1}+\frac{\theta}{\lambda-t}.$$ If now $\lambda$ and $\mu$ are viewed as functions of $t$, substituting this condition into the Hamiltonian system leading to the sixth Painlev\'e equation, we get that $\lambda$ satisfies the Riccati equation $$\frac{d\lambda}{dt}=\frac{\kappa_0 t-(1+\kappa_0+(\kappa_0+\kappa_1)t+\theta)\lambda+(1+\kappa_0+\kappa_1+\theta)\lambda^2}{t(t-1)}$$ and $\kappa_0+\kappa_1+\theta+\kappa=0$. This gives classical solutions of the sixth Painlev\'e equation provided that $\kappa_0=\pm \kappa_{\infty}-\theta-\kappa_1-1$. However, with this additional condition on the parameters we have $\alpha\beta=0$ and $q=0$.

In the equation for the derivative of the confluent Heun function (\ref{dCHE}) we first make the change of variables $v(z)\to (1-z/(z-1))^{\sigma}v(z/(z-1))$ and renaming the new independent variable as $z$ again, we put
\begin{gather*}
\gamma=-\kappa_0,\;\;\delta=\kappa_0+\theta+2\sigma,\;\;\varepsilon=-t\eta,\\
\sigma=-\frac{1}{2}(\kappa_0\pm \kappa_{\infty}+\theta),\;\;q=\frac{\alpha \lambda}{\lambda-1},\;\;\alpha=\frac{1}{2}t\eta(2+\kappa_0\pm\kappa_{\infty}+\theta).
\end{gather*}
The resulting equation  is the same as equation (\ref{eq}) with (\ref{P5 coeff p1}), (\ref{P5 coeff p2}) and the expression for $H_{V}$ provided that
$$\mu=\frac{\kappa_0}{\lambda}-\frac{t\eta}{(\lambda-1)^2}+\frac{\theta-\kappa_0\pm \kappa_{\infty}}{2(\lambda-1)}.$$
Substituting this condition into the Hamiltonian system leading to the fifth Painlev\'e equation, we get that $\lambda$ satisfies the Riccati equation $$t\frac{d\lambda}{dt}\pm \kappa_{\infty}\lambda^2-(\pm \kappa_{\infty}-\kappa_0-t\eta)-\kappa_0=0$$ and $\eta(2+\kappa_0\pm\kappa_{\infty}+\theta)=0$.  Again, with this additional condition on the parameters we have $\alpha=0$ and $q=0$.

In the equation for the derivative of the  bi-confluent Heun function  (\ref{dBHE1}) we take
\begin{gather*}
\gamma=-\kappa_0,\;\;\delta=-t,\;\;q=\alpha\lambda,\;\;
\alpha=\frac{\theta_{\infty}+1}{2},\;\;\varepsilon=-\frac{1}{2}.
\end{gather*}
The resulting equation  is the same as equation (\ref{eq}) with (\ref{P4 coeff p1}), (\ref{P4 coeff p2}) and the expression for $H_{IV}$ provided that
$$\mu=t+\frac{\kappa_0}{\lambda}+\frac{\lambda}{2}.$$
Substituting this condition into the Hamiltonian system leading to the fourth Painlev\'e equation, we get that $\lambda$ satisfies the Riccati equation $$\frac{d\lambda}{dt}=\lambda^2+2t\lambda+2\kappa_0$$ and  $\theta_{\infty}+1=0$.  Again, with this additional condition on the parameters we have $\alpha=0$ and $q=0$.

In the equation for the derivative of the  double-confluent Heun function  (\ref{dBHE1}) we take
\begin{gather*}
\gamma=t\eta_0,\;\;\delta=-1-\theta_0,\;\;q=\alpha\lambda,\;\;
\alpha=\frac{1}{2}\eta_{\infty}(\theta_0+\theta_{\infty}+2),\;\;\varepsilon=-\eta_{\infty}.
\end{gather*}
The resulting equation  is the same as equation (\ref{eq}) with (\ref{P3 coeff p1}), (\ref{P3 coeff p2}) and the expression for $H^{\prime}_{III}$ provided that
$$\mu=\eta_{\infty}-\frac{t\eta_0}{\lambda^2}+\frac{\theta_0+1}{\lambda}.$$
Substituting this condition into the Hamiltonian system leading to the modified third Painlev\'e equation  $P^{\prime}_{III}$, we get that $\lambda$ satisfies the Riccati equation $$t\frac{d\lambda}{dt}= \eta_{\infty}\lambda^2+(\theta_0+2)\lambda-t\eta_0$$ and $\eta_{\infty}(\theta_0+\theta_{\infty}+2)=0$.  Again, with this additional condition on the parameters we have $\alpha=0$ and $q=0$.

In the equation for the derivative of the  tri-confluent Heun function  (\ref{dTHE}) we take
\begin{gather*}
\gamma=-t,\;\;\delta=0,\;\;q=\alpha\lambda,\;\;
\alpha=1-2\alpha_2,\;\;\varepsilon=-2.
\end{gather*}
The resulting equation  is the same as equation (\ref{eq}) with (\ref{P2 coeff p1}), (\ref{P2 coeff p2}) and the expression for $H_{II}$ provided that
$$\mu=2\lambda^2+t.$$
Substituting this condition into the Hamiltonian system leading to the second Painlev\'e equation, we get that $\lambda$ satisfies the Riccati equation $$2\frac{d\lambda}{dt}=2 \lambda^2+t$$ and $2\alpha_2=1$.  Again, with this additional condition on the parameters we have $\alpha=0$ and $q=0$.

Hence, we see that in all cases we can reduce equations for the derivatives of the Heun functions to linear equations isomonodromy deformations of which lead to the Painlev\'e equations with an additional constraint on $\lambda $ and $\mu$. However, in order to get classical solutions of the Painlev\'e equations we need an additional constraint on the parameters. Therefore, those linear equations isomonodromy deformations of which are described by classical solutions of the Painlve\'e equations cannot be obtained from the equations for the derivatives of the Heun functions.

\section*{Acknowledgements}

We thank M. Nieszporski (University of Warsaw)  for interesting discussions. GF acknowledges the support of the National Science Center (Poland) via grant OPUS 2017/25/ B/BST1/00931 and the Alexander von Humboldt Foundation. The support of the Armenian State Committee of Science (SCS Grants No. 18RF-139 and No. 18T-1C276), the Armenian National Science and Education Fund (ANSEF Grant No. PS-4986), the Russian-Armenian (Slavonic) University is also greatfully acknowledged. AI thanks the colleagues from the University of Warsaw for hospitality and inspiring discussions.

\end{document}